\documentclass[a4paper,10pt]{amsart}
\usepackage[utf8]{inputenc}

\newtheorem{thm}{Theorem}[section]
\newtheorem*{thm*}{Theorem}

\newtheorem{cor}[thm]{Corollary}
\theoremstyle{definition}
\newtheorem{definition}[thm]{Definition}

\newtheorem{remark}[thm]{Remark}
\theoremstyle{remark}
\newcommand{\bbZ}{{{\mathbb{Z}}}}

\newcommand{\bbP}{{{\mathbb{P}}}}

\newcommand{\vg}{{{\lbrace{0,1\rbrace}^{V(G)}}}}

\title{Mixing of The Noisy Voter Model}
\author{Harishchandra Ramadas}
\thanks{University of Washington, Seattle, USA. Contact: \texttt{ramadas.math.washington.edu}}

\begin{document}

\begin{abstract}
We prove that the noisy voter model mixes extremely fast -- in time of $O(\log n)$ on any graph with $n$ vertices -- for arbitrarily small values of the `noise parameter'. We then explain why, as a result, this is an example of a spin system that is \emph{always} in the `high-temperature regime'.
\end{abstract}

\maketitle

\section{Introduction}\label{Intro}

The \emph{noisy voter model} is a spin system that can serve as a crude stochastic model for the spread of opinions in a human population, or the competition for territory between two species. The model was introduced by Granovsky and Madras \cite{GM} as a simple ergodic variant of the extensively-studied \emph{voter model} \cite{Lig0}. Our main result in this paper (Theorem \ref{nvmotm}) is an upper bound on the mixing time of this model: we show that on an arbitrary graph with $n$ vertices, the system mixes extremely fast -- in time of $O(\log n)$ -- for \emph{all} values of the `noise parameter', $\delta$.

It is well-known that fast \emph{temporal} mixing is equivalent to a \emph{spatially} `well-mixed' stationary measure \cite{DSVW}. Corollary \ref{ssmnvm} states that the stationary measure of the noisy voter model on any subset of the integer lattice satisfies a decay-of-correlations condition called \emph{strong spatial mixing}. Conventional spin models such as the Ising model, gas hard-core model and proper coloring model exhibit strong spatial mixing only when they are sufficiently `disordered', e.g. the Ising model on boxes in $\bbZ^2$ above the critical temperature. The noisy voter model is thus extremely special, since it can be thought of as \emph{always} being in the high-temperature regime.

\section{Preliminaries}

\subsection{The Noisy Voter Model}\label{SSVM}

Starting with a finite simple undirected graph $G$ with vertex set $V(G)$, the state space we consider is the set of all $0$-$1$ configurations on $V(G)$, denoted by $\lbrace{0,1\rbrace}^{V(G)}$. In other words, at each site $v \in V(G)$, we place a `spin', which may have value $0$ or $1$. For $\eta \in \vg$ and $v \in V(G)$, we write $\eta(v)$ for the value of the spin at the site $v$. For $x \in V(G)$, we denote by $\eta^x$ the configuration obtained from $\eta$ by flipping the spin at $x$, leaving all other spins unchanged:

\begin{equation}\label{**}
\eta^x(v) :=
    \begin{cases}
      1-\eta(v) & v=x,\\
      \eta(v) & v \neq x.\\
    \end{cases}
\end{equation}

Further, we denote by $d(x)$ the degree of the vertex $x$ and write $y \sim x$ to indicate that the vertices $x$ and $y$ are neighbours. 

\begin{definition}
 Given a finite graph $G$ and a function
 \[c: V(G)\times\vg \to [0,\infty),\]
  a \emph{spin system} is a continuous-time Markov chain on the state space $\vg$ with transition-rate matrix $q(\cdot,\cdot)$, where
\begin{equation*}
q(\eta,\eta')=
    \begin{cases}
      c(x,\eta) & $if $\eta' = \eta^x $ for some $ x\in V(G),\\
      0 & $otherwise,$\\
    \end{cases}
\end{equation*}
for all configurations $\eta, \eta' \in \vg$ such that $\eta \neq \eta'$.
\end{definition}

The \emph{voter model}, introduced in the 70s, is an extremely well-understood spin system \cite[Chapter V]{Lig0}. It is defined by the transition rates
\[c(x,\eta) = \frac{1}{d(x)}\#\lbrace{y:y\sim x, \eta(y)\not=\eta(x)\rbrace},\]
i.e., the rate of flipping the spin at $x$ while in state $\eta$ is proportional to the number of neighbours of $x$ with spin differing from that at $x$.

We imagine a `voter' located at each vertex $x$ and that the spin at $x$ corresponds to its `opinion', either $0$ or $1$. Another interpretation of the voter model is a situation where two species, denoted by $0$ and $1$, compete for territory.

Intuitively, we imagine independent rate-$1$ Poisson clocks located at each vertex $v \in V(G)$ (i.e., clocks that ring after independent $\mathrm{Exp}(1)$ waiting times). If the clock at vertex $x$ rings, the voter at $x$ picks one of its neighbours uniformly at random, and modifies its own opinion to match that of the chosen neighbour.

We observe, however, that the voter model is \emph{not} ergodic, since the all-$0$ and all-$1$ states (denoted $\underbar{0}$ and $\underbar{1}$ respectively) are absorbing. There are hence multiple invariant measures: $\delta_{\underbar{0}}$, $\delta_{\underbar{1}}$ and all their convex combinations.

The model we study in this paper is a simple ergodic variant of the voter model -- the \textbf{noisy voter model}, introduced by Granovsky and Madras in 1995 \cite{GM}.

\begin{definition}
Let $\delta > 0$. The \emph{noisy voter model} is a spin system on a graph $G$ with transition rates
\begin{equation}\label{1}
c(x,\eta)=\frac{1}{2\delta+1}\left[\frac{1}{d(x)}\#\{y:y\sim x, \eta(y)\not=\eta(x)\}+\delta\right],
\end{equation}
where $d(x)$ is the degree of the vertex $x$.
\end{definition}

The factor $\frac{1}{2\delta+1}$ in \eqref{1} `normalises' the system so that we have the following intuitive description of the dynamics: assign independent rate-$1$ Poisson clocks to each vertex $v \in V(G)$. Every time a clock rings at a vertex $v$, erase $\eta(v)$ and set

\begin{itemize}
\item $\eta(v) = 0$ with probability $\frac{1}{2\delta+1}\left[\frac{1}{d(v)}\#\{y:y\sim x , \eta(y)=0\} + \delta\right]$,
\item $\eta(v) = 1$ with probability $\frac{1}{2\delta+1}\left[\frac{1}{d(v)}\#\{y:y\sim x , \eta(y)=1\} + \delta\right]$.
\end{itemize}

The parameter $\delta$ is called the `noise' for obvious reasons, and we assume that it is strictly positive. (The case $\delta = 0$ is just the regular voter model, which, as explained above, is not ergodic.)

\begin{remark}
 We mention briefly that this is not the most general form of the noisy voter model -- we could alternatively have \emph{two} positive parameters $\delta$ and $\beta$ so that the noise in the model causes the system to preferentially tend towards spin $0$ or spin $1$. In this case, the flip rates can be written as follows:
\[c(x,\eta)=\frac{1}{\delta  + \beta +1}\left[\frac{1}{d(x)}\#\{y:y\sim x, \eta(y)\not=\eta(x)\}+\delta\textbf{1}_{\eta(x)=0} + \beta\textbf{1}_{\eta(x)=1}\right].\]
However, for the sake of simplicity of the proof, we choose to restrict ourselves to a single noise parameter, $\delta$; it can be easily verified that the proof goes through even in the more general case.
\end{remark}

\subsection{Mixing Time}

We first formally define the intuitive notion of `distance to equilibrium' for an ergodic finite-state Markov chain, or more generally, the distance between two probability measures on a finite set.

\begin{definition}
For any two probability measures $\mu_1$ and $\mu_2$ on a finite set $S$, the \emph{total variation distance} between $\mu_1$ and $\mu_2$ is defined to be
\begin{equation}\label{TV}
 \|\mu_1 - \mu_2\| := \max_{A \subseteq S} |\mu_1(A) - \mu_2(A)|.
\end{equation}
\end{definition}

We define next the notion of the mixing time of a Markov Chain (as a function of $\epsilon$), the minimum time it takes to get to within a total variation distance of $\epsilon$ to stationarity:

\begin{definition}
For $0 < \epsilon < 1$, the (total variation) \emph{mixing time} of an ergodic continuous-time Markov chain $(X_t)_{t \geq 0}$ on a finite state space $S$ with stationary distribution $\mu$ is defined to be
\[
 t_{\mathrm{mix}}(\epsilon) := \inf \lbrace{t \geq 0: \max_{x \in S} \|\mathbb P_x(X_t \in \cdot) - \mu\| \leq \epsilon \rbrace}.
\]
\end{definition}

For two random variables $X$ and $Y$ with values in the same state space, we write $\|X - Y\|$ for the total variation distance between their distributions.

\begin{definition}\label{OTMdef}
 We say that a (continuous-time) spin system has \emph{optimal temporal mixing} if there exist constants $C$, $c > 0$ such that the dynamics on any finite graph $G$ has the following property: for any two copies $(X_{t})_{t\geq0}$ and $(Y_{t})_{t\geq0}$ of the chain with possibly different initial configurations, we have	
\begin{equation}\label{2}
\|X_{t} - Y_{t}\| \leq Cne^{-ct},
\end{equation}
where, as mentioned above, the left hand side refers to the total variation distance between the distributions of the $\lbrace0, 1\rbrace^{V(G)}$-valued random variables $X_t$ and $Y_t$, and $n = |V(G)|$, the number of vertices in G.
\end{definition}

\begin{remark}\label{otmlogn}
We observe that optimal temporal mixing is equivalent to a mixing time of $O\left(\log \left(\frac{n}{\epsilon}\right)\right)$:
\begin{itemize}
 \item Assume first that optimal temporal mixing holds, i.e., for two instances $(X_t)_{t\geq0}$ and $(Y_t)_{t\geq0}$ of the dynamics, $\|X_t - Y_t\| \leq Cne^{-ct}$. We further assume that $Y_0 {\buildrel d \over =} \mu$, where $\mu$ is the equilibrium measure of the dynamics. Hence $Y_t {\buildrel d \over =} \mu$ for all $t \geq 0$. We thus obtain $\|X_t - \mu\| \leq Cne^{-ct}$. Now, suppose $t = \frac{1}{c}\log\left(\frac{nC}{\epsilon}\right) = O \left(\log\left(\frac{n}{\epsilon}\right)\right)$. Then, clearly $\|X_t - \mu\| \leq \epsilon$.
 \item Conversely, assume that $t_{\mathrm{mix}}(\epsilon) \leq C \log\left(\frac{n}{\epsilon}\right)$,
and let $(X_t)_{t\geq0}$ and $(Y_t)_{t\geq0}$ be two instances of the dynamics. Then, for all $t > 0$, $t = C \log \left(\frac{2n}{\epsilon}\right)$ implies $\|X_t - \mu\| \leq \frac{\epsilon}{2}$. But $t = C \log \left(\frac{2n}{\epsilon}\right) \Leftrightarrow \epsilon = 2ne^{-t/C}$. Hence, for all $t > 0$, 
\[\|X_t - Y_t\| \leq \|X_t - \mu\| + \|Y_t - \mu\| < \frac{\epsilon}{2} + \frac{\epsilon}{2} = \epsilon = 2ne^{-t/C},\] i.e., optimal temporal mixing holds.
\end{itemize}
\end{remark}

Throughout the rest of this article, we keep the parameter $\epsilon$ fixed.

\section{Main Result: Optimal Temporal Mixing}

\begin{thm}\label{nvmotm}
 For any $\delta > 0$, the noisy voter model on any finite graph $G$ with $n$ vertices has mixing time of $O(\log n)$.

\begin{proof}

A standard method for obtaining an upper bound on the total variation distance between the distributions of two random variables $X$ and $Y$ is to construct a suitable coupling $(X, Y)$, and to use the well-known inequality,  
\begin{equation}\label{3}
\|X - Y\| \leq \mathbb P(X \neq Y),
\end{equation}
which holds for \emph{any} coupling $(X, Y)$ of $X$ and $Y$ (see, e.g. \cite[Proposition 4.7]{LPW}).

Further, if $(X_t, Y_t)_{t \geq 0}$ is a coupled spin system defined on a finite graph $G$, then for a fixed time $t$, a simple union bound yields
\begin{equation}\label{4}
\begin{aligned}
\mathbb P(X_t \neq Y_t) &= \mathbb P\left(\bigcup_{v\in V(G)}\left\lbrace X_t(v) \neq Y_t(v)\right\rbrace\right) \\
&\leq \sum_{v\in V(G)}\mathbb P(X_t(v) \neq Y_t(v)) \\
&\leq n\left(\max_{v\in V(G)}\mathbb P\left(X_t(v) \neq Y_t(v)\right)\right).
\end{aligned}
\end{equation}
Hence, from \eqref{2}, \eqref{3} and \eqref{4}, it is enough to show that there exists a coupling $(X_t, Y_t)_{t\geq0}$ such that for arbitrary starting configurations $X_0$ and $Y_0$, the  following holds for all $t \geq 0$:
\begin{equation}\label{5}
\max_{v\in V(G)}\mathbb P\left( X_t(v) \neq Y_t(v)\right) \leq Ce^{-ct},
\end{equation}
where $C$ and $c$ are fixed constants independent of $G$. 

\textbf{The Coupling:} Since the noisy voter model is a monotone system, it suffices to consider a coupling $(X_t, Y_t)$ of the system, with $X_0 = \underbar{0}$ and $Y_0 = \underbar{1}$. Place independent rate-1 Poisson clocks at each vertex $v \in V(G)$. When the clock at $v$ rings:
\begin{itemize}
 \item With probability $\frac{\delta}{2\delta+1}$, set both $X_t(v), Y_t(v) = 0$.
 \item With probability $\frac{\delta}{2\delta+1}$, set both $X_t(v), Y_t(v) = 1$.
 \item With probability $\frac{1}{2\delta+1}$, pick a neighbour $u$ of $v$ uniformly at random. Set $X_t(v) = X_t(u)$ and $Y_t(v) = Y_t(u)$.
\end{itemize}
It is easy to see that both $X_t$ and $Y_t$ follow the dynamics of the noisy voter model, defined in \eqref{1}.

Let $Z_t(v) := \mathbf{1}_{\lbrace X_t(v) \neq Y_t(v) \rbrace}$ be the \emph{disagreement process}, i.e., the indicator of the event that the processes starting from the $\underbar{0}$ and $\underbar{1}$ configurations ($X_t$ and $Y_t$ respectively) disagree at vertex $v$ at some time $t\geq0$. Obviously, $Z_0(v) = 1$  for all $v \in V(G)$.

We hence wish to show that
\begin{equation}\label{6}
\max_{v\in V(G)}\mathbb P\left(Z_t(v) = 1\right) \leq Ce^{-ct}.
\end{equation}

With this coupling, every time the Poisson clock at vertex $v$ rings, with probability at least $\frac{2\delta}{2\delta+1}$, both $X_t$ and $Y_t$ are updated at $v$ \emph{with the same value} ($0$ or $1$). Hence, independent of its current value, $Z_t(v)$ is set to $0$. With probability $\frac{1}{2\delta+1}$, the vertex $v$ chooses a neighbour $u$ uniformly at random. Independent of the current value of $Z_t(v)$, it is clear that if $Z_t(u) = 1$, then $Z_t(v)$ becomes $1$ whereas if $Z_t(u) = 0$, then $Z_t(v)$ becomes $0$.

Since $V(G)$ is a finite set, $(Z_t)_{t\geq0}$ is just a continuous-time Markov chain on the (finite) state space $\lbrace0,1\rbrace^{V(G)}$ with initial state $\underbar{1}$. We denote the corresponding transition function $p_t(\cdot,\cdot)$ and transition-rate matrix $q(\cdot,\cdot)$. As before, for $\eta \in \lbrace0,1\rbrace^{V(G)}$, denote by $\eta^{v}$ the configuration with the value at $v$ flipped, and everywhere else the same as $\eta$. We can hence write a formula for the transition rates similar to the noisy voter model \eqref{1}:
\[ q(\eta,\eta^{v}) = \frac{2\delta}{2\delta+1}\textbf{1}_{\lbrace\eta(v)=1\rbrace} + \frac{1}{2\delta+1}\cdot\frac{1}{d(v)}\sum_{u: u\sim v} \textbf{1}_{\lbrace\eta(u)\neq\eta(v)\rbrace}.\]
Interchanging $\eta$ and $\eta^{v}$, we get
\[ q(\eta^{v},\eta) = \frac{2\delta}{2\delta+1}\textbf{1}_{\lbrace\eta(v)=0\rbrace} + \frac{1}{2\delta+1}\cdot\frac{1}{d(v)}\sum_{u: u\sim v} \textbf{1}_{\lbrace\eta(u)=\eta(v)\rbrace}.\]
Since $q(\eta,\eta')=0$ if $\eta$ and $\eta'$ differ at more than one site and since for any Markov chain $\sum_{y}q(x,y)=0$, we have
\begin{align*}
 q(\eta,\eta) &= - \sum _{v\in V(G)} q(\eta,\eta^{v}) \\
&= - \sum _{v\in V(G)} \left(\frac{2\delta}{2\delta+1}\textbf{1}_{\eta(v)=1} + \frac{1}{2\delta+1}\cdot\frac{1}{d(v)}\sum_{u: u\sim v} \textbf{1}_{\lbrace\eta(u)\neq\eta(v)\rbrace}\right).
\end{align*}
We can now write the forward Kolmogorov equation (see, e.g. \cite[Equation 2.17]{Lig}) for the Markov chain $(Z_t)$, started at $\underbar{1}$:
\small
\begin{equation}\label{7}
\begin{split}
&\frac{d}{dt}p_{t}(\underbar{1},\eta) = \sum_{v \in V(G)} \left[p_{t}(\underbar{1},\eta^{v}) \cdot q(\eta^{v},\eta)\right] + p_{t}(\underbar{1},\eta) \cdot q(\eta,\eta) \\
&= \frac{2\delta}{2\delta+1}\sum_{v \in V(G)} p_{t}(\underbar{1},\eta^{v})\textbf{1}_{\lbrace\eta(v)=0\rbrace} 
+ \frac{1}{2\delta+1}\sum_{v \in V(G)} \left( \frac{1}{d(v)} \cdot p_{t}(\underbar{1},\eta^{v}) \sum_{u: u \sim v} \textbf{1}_{\lbrace\eta(u)=\eta(v)\rbrace} \right) \\
&- \frac{2\delta}{2\delta+1}\sum_{v \in V(G)} p_{t}(\underbar{1},\eta)\textbf{1}_{\lbrace\eta(v)=1\rbrace} 
- \frac{1}{2\delta+1}\sum_{v \in V(G)} \left( \frac{1}{d(v)} \cdot p_{t}(\underbar{1},\eta) \sum_{u: u \sim v} \textbf{1}_{\lbrace\eta(u) \neq \eta(v)\rbrace} \right).
\end{split}
\end{equation}
\normalsize
Now, for $x \in V(G)$,
\begin{equation}\label{7.1}
\frac{d}{dt}\mathbb P(Z_t(x)=1) = \frac{d}{dt}\sum_{\eta \in \lbrace0,1\rbrace^{V(G)}} p_{t}(\underbar{1},\eta) \cdot \textbf{1}_{\lbrace\eta(x)=1\rbrace}. 
\end{equation}
Substituting \eqref{7} in \eqref{7.1} above, with the convention that sums indexed with $v$ and $\eta$ run over all vertices and all configurations respectively,
\small
\begin{equation}\label{7.2}
\begin{aligned}
 &\frac{d}{dt}\bbP(Z_{t}(x)=1) \\
 &=\frac{2\delta}{2\delta+1} \sum_{v,\eta}\left(p_{t}(\underbar{1},\eta^{v})\textbf{1}_{\lbrace\eta(v)=0,\eta(x)=1\rbrace} - p_{t}(\underbar{1},\eta)\textbf{1}_{\lbrace\eta(v)=1,\eta(x)=1\rbrace} \right) \\
 &+ \frac{1}{2\delta+1} \sum_{v,\eta} \frac{1}{d(v)} \cdot \textbf{1}_{\lbrace\eta(x)=1\rbrace} \cdot
 \left(p_{t}(\underbar{1},\eta^{v}) \sum_{u: u\sim v} \textbf{1}_{\lbrace\eta(u)=\eta(v)\rbrace}
 - p_{t}(\underbar{1},\eta) \sum_{u: u\sim v} \textbf{1}_{\lbrace\eta(u)\neq\eta(v)\rbrace} \right) .
\end{aligned}
\end{equation}
\normalsize
For each of the sums over $v$ in \eqref{7.2} above, we consider the cases $v=x$ and $v \neq x$ separately to get
\small
\begin{align*}
 &\frac{d}{dt}\mathbb P(Z_{t}(x)=1) \\
 &=\frac{2\delta}{2\delta+1} \sum_{\eta}\left(p_{t}(\underbar{1},\eta^{x}) \cdot 0 - p_{t}(\underbar{1},\eta)\cdot\textbf{1}_{\lbrace\eta(x)=1\rbrace} \right) \\
 &+\frac{2\delta}{2\delta+1} \sum_{v: v\neq x} \sum_{\eta}\left(p_{t}(\underbar{1},\eta^{v})\cdot \textbf{1}_{\lbrace\eta^{v}(v)=1,\eta^{v}(x)=1\rbrace} - p_{t}(\underbar{1},\eta)\cdot \textbf{1}_{\lbrace\eta(v)=1,\eta(x)=1\rbrace} \right) \\
 &+ \frac{1}{2\delta+1} \bigg\{\sum_{\eta} \frac{1}{d(x)} \cdot \textbf{1}_{\lbrace\eta(x)=1\rbrace} \cdot
 \left(p_{t}(\underbar{1},\eta^{x}) \sum_{u: u\sim x} \textbf{1}_{\lbrace\eta(u)=1\rbrace}
 - p_{t}(\underbar{1},\eta) \sum_{u: u\sim x} \textbf{1}_{\lbrace\eta(u)=0\rbrace} \right) \\
 &+ \sum_{v: v\neq x} \frac{1}{d(v)} \sum_{\eta} \textbf{1}_{\lbrace\eta(x)=1\rbrace} \cdot
 \left(p_{t}(\underbar{1},\eta^{v}) \sum_{u: u\sim v} \textbf{1}_{\lbrace\eta(u)=\eta(v)\rbrace}
 - p_{t}(\underbar{1},\eta) \sum_{u: u\sim v} \textbf{1}_{\lbrace\eta(u)\neq\eta(v)\rbrace} \right) \bigg\}.
\end{align*}
\normalsize
In the second term on the right hand side, the sum over $\eta$ is clearly zero. Since in the second term in the braces $\lbrace{\ldots\rbrace}$, $v \neq x$, it can be rewritten as
\begin{align*}
\sum_{v: v\neq x} \frac{1}{d(v)} \sum_{\eta} \bigg(&p_{t}(\underbar{1},\eta^{v}) \cdot \textbf{1}_{\lbrace\eta^v(x)=1\rbrace} \sum_{u: u\sim v} \textbf{1}_{\lbrace\eta^{v}(u)\neq\eta^{v}(v)\rbrace} \\
 - &p_{t}(\underbar{1},\eta) \cdot \textbf{1}_{\lbrace\eta(x)=1\rbrace}\sum_{u: u\sim v} \textbf{1}_{\lbrace\eta(u)\neq\eta(v)\rbrace} \bigg),
\end{align*}
and as before, the sum over $\eta$ equals zero. We are hence left with
\small
\begin{align*}
 &\frac{d}{dt}\mathbb P(Z_{t}(x)=1) \\
 &=\frac{2\delta}{2\delta+1} \sum_{\eta}\left(-p_{t}(\underbar{1},\eta)\cdot\textbf{1}_{\lbrace\eta(x)=1\rbrace} \right) \\
 &+ \frac{1}{2\delta+1} \sum_{\eta} \frac{1}{d(x)} \cdot \textbf{1}_{\lbrace\eta(x)=1\rbrace} \cdot
 \left(p_{t}(\underbar{1},\eta^{x}) \sum_{u: u\sim x} \textbf{1}_{\lbrace\eta(u)=1\rbrace}
 - p_{t}(\underbar{1},\eta) \sum_{u: u\sim x} \textbf{1}_{\lbrace\eta(u)=0\rbrace} \right) .
 \end{align*}
\normalsize
We finally write the transition function elements and indicator functions in terms of probabilities to get the following `master equation', valid for all $x \in V(G)$:
\begin{equation}
\begin{aligned}\label{8}
\frac{d}{dt}\mathbb P(Z_t(x)=1) &= \frac{-2\delta}{2\delta+1}\mathbb P(Z_t(x)=1) \\
&- \frac{1}{2\delta+1}\cdot\frac{1}{d(x)}\left(\sum_{u: u \sim x} \mathbb P(Z_t(x)=1, Z_t(u)=0)\right) \\
&+ \frac{1}{2\delta+1}\cdot\frac{1}{d(x)}\left(\sum_{u: u \sim x} \mathbb P(Z_t(x)=0, Z_t(u)=1)\right).
\end{aligned}
\end{equation}
For notational simplicity, we now define
\begin{equation}\label{Mt}
M(t) := \max_{1\leq k\leq n} \bbP\left(Z_t(x_k) = 1\right),
\end{equation}
where $x_1,\ldots,x_n$ are the $n$ vertices of $G$.
Suppose, at some instant $s\geq 0$,
\[M(s) = \bbP(Z_s(x_k)=1),\]
for some $k$ with $1\leq k \leq n$. Then $\bbP\left(Z_s(x_k)=1\right) \geq \bbP\left(Z_s(x_l)=1\right)$ for all $l \neq k$, $1 \leq l \leq n$. Hence,
\begin{align*}
\bbP(Z_s(x_k)=1, Z_s(x_l)=0) &= \bbP(Z_s(x_k)=1) - \bbP(Z_s(x_k)=1, Z_s(x_l)=1) \\
&\geq \bbP(Z_s(x_l)=1) - \bbP(Z_s(x_k)=1, Z_s(x_l)=1) \\
&= \bbP(Z_s(x_k)=0, Z_s(x_l)=1).
\end{align*}
Setting $t=s$ and $x=x_k$ in \eqref{8}, we see that the magnitude of the second term on the right hand side is greater than that of the third term, since for all neighbours $u$ of $x_k$,
\[\bbP(Z_s(x_k)=1, Z_s(u)=0) \geq \bbP(Z_s(x_k)=0, Z_s(u)=1).\]
As a result,
\begin{equation}\label{MsIn}
\begin{aligned}
\frac{d}{dt}\bbP(Z_t(x_k)=1) \Bigr\rvert_{t=s} &\leq -\frac{2\delta}{2\delta+1} \bbP(Z_s(x_k)=1)\\
&= -\frac{2\delta}{2\delta+1}M(s).
\end{aligned}
\end{equation}
In general, at the instant $t=s$, there might be multiple $k$ ($1\leq k\leq n$), such that $M(s) = \bbP(Z_s(x_k) = 1)$. However, \eqref{MsIn} holds for \emph{every} such $k$. Furthermore, $\bbP(Z_t(x_l)=1)$ is a differentiable (hence continuous) function of $t$ for each $l=1,\ldots,n$. Consequently, for all $t \geq 0$, we have the inequality 
\[D^+M(t) := \limsup_{h \to 0^+} \frac{M(t+h)-M(t)}{h} \leq -\frac{2\delta}{2\delta+1}M(t),\]
where $D^+M(t)$ is the \emph{upper right Dini derivative} of $M$ at $t$. Again, since $\bbP(Z_t(x_l)=1)$ is continuous for each $l$, $M(t)$ is continuous at every $t \geq 0$, and we can use the fundamental theorem of calculus to write
\[M(t) = \frac{d}{dt}\int_0^t M(s) ds = D^+ \left(\int_0^t M(s)\ ds\right).\]
Hence,
\[D^+M(t) + \frac{2\delta}{2\delta+1}D^+ \left(\int_0^t M(s)\ ds\right) \leq 0.\]
Using the elementary fact that $D^+(f+g) \leq D^+f + D^+g$ for any two functions $f$ and $g$, we have
\[D^+\left(M(t) + \frac{2\delta}{2\delta+1} \int_0^t M(s)\ ds\right) \leq 0,\]
for all $t \geq 0$. 
The monotonicity theorem for Dini derivatives \cite[Appendix I]{RHL} says that if $D^+f(t) \leq 0$ for all $t \geq 0$, $f$ must be nonincreasing. As a result,
\[M(t) + \frac{2\delta}{2\delta+1} \int_0^t M(s)\ ds\leq M(0) = 1,\]
and hence,
\[M(t) \leq 1 - \frac{2\delta}{2\delta+1} \int_0^t M(s)\ ds\]
for all $t \geq 0$. Finally, applying the Gronwall-Bellman inequality \cite{B}, we get
\[M(t) \leq e^{-\frac{2\delta}{2\delta+1}t},\]
i.e.,
\[\max_{v\in V(G)} \bbP(Z_t(v) = 1) \leq e^{-\frac{2\delta}{2\delta+1}t}.\]
We thus conclude that optimal temporal mixing holds, or equivalently, by Remark \ref{otmlogn}, that the system mixes in time of $O(\log n)$. 
\end{proof}
\end{thm}

\begin{remark}\label{NVMtheta}
It is generally believed that for any spin system, $O(\log n)$ mixing is indeed `optimal': by a coupon-collector argument, one imagines that in $o(\log n)$ time, not enough sites are updated for the system to equilibriate. However, it is shown in \cite{HS} that this argument can fail, but the authors then prove that any \emph{reversible} spin system with `local' interactions on a bounded-degree graph with $n$ vertices necessarily has mixing time of $\Omega(\log n)$. Unfortunately, a simple application of Kolmogorov's criterion \cite[Chapter 1]{K} shows that the noisy voter model is in general \emph{not} reversible. An exception is the simple case of the 1-dimensional torus, $\bbZ/n\bbZ$, for which it is easy to see that the noisy voter model with parameter $\delta$ is in fact exactly the same process as the stochastic (ferromagnetic) Ising model with inverse temperature $\beta = \frac14\log(1+\delta^{-1})$.
\end{remark}

\section{Always Hot: Strong Spatial Mixing}
As mentioned in the introduction, it is already well-known that fast temporal mixing is in some sense equivalent to a `well-mixed' stationary measure. We formalise this below.

Although our main result (optimal temporal mixing) holds for any graph on $n$ vertices, we focus on the special case of `boxes' in the $d$-dimensional integer lattice $\bbZ^d$, for some $d \geq 1$. (We do this for simplicity; in general, one could consider any lattice with subexponential growth.) By a (finite) \emph{box} in $\bbZ^d$, we mean a finite induced subgraph: a finite subset of the vertex set, along with all edges between these vertices inherited from the lattice. (In particular, a box need not be cuboidal.) 

For $u, v \in \mathbb Z^d$, the \emph{distance} between $u$ and $v$ is defined to be
\[ \mathrm{dist}(u,v) := \sum_{i=1}^{d} |u_i - v_i|.\]

The distance between a vertex $u$ and a box $G$ is
 \[\mathrm{dist}(u,G) := \min_{v\in V(G)}\mathrm{dist}(u,v).\]

For a box $G \subseteq \bbZ^d$, we define its \emph{boundary} to be the set of vertices not in $G$ that are neighbours with at least one vertex in $G$. A configuration $\tau \in \lbrace{0,1\rbrace}^{\partial G}$ is said to be a \emph{boundary condition} for the spin system dynamics on the box $G$, if we run the dynamics on $G \cup \partial G$ conditioned on the spins at the boundary being frozen at their initial values. 

Finally, for a finite box $H \subseteq G \subset \mathbb Z^d$ and a boundary condition $\tau \in \lbrace{0,1\rbrace}^{\partial G}$, we define $\mu_{G}^{\tau}|_H$ to be the stationary distribution for the dynamics run on $G$ with boundary condition $\tau$, \emph{projected onto} $\lbrace{0,1\rbrace}^{V(H)}$.

Dyer, Sinclair, Vigoda and Weitz give in \cite{DSVW} an elegant combinatorial proof of the fact that if a spin system has optimal temporal mixing, then the stationary distribution satisfies a condition called strong spatial mixing, which we define below.

\begin{definition}\label{ssm}
 A spin system has strong spatial mixing if there exist constants $C, c > 0$ such that for any  finite box $G \subseteq \bbZ^d$ and a box $H \subseteq G$, any site $u \in \partial G$ and any pair of boundary conditions $\tau$ and $\tau^u$ that differ only at $u$,
\begin{equation}\label{ssmdef}
 \| \mu_{G}^{\tau}|_H - \mu_{G}^{\tau^u}|_H \| \leq C|V(H)|\exp\left(-c\cdot \mathrm{dist}(u,H)\right)
\end{equation}
 holds.
\end{definition}

Due to \cite[Theorem 2.3]{DSVW}, we hence immediately have the following corollary of Theorem \ref{nvmotm}:

\begin{cor}\label{ssmnvm}
The noisy voter model has strong spatial mixing for \emph{all} values of the noise parameter $\delta$.
\end{cor}

For other spin systems like the Ising model (with zero magnetic field), proper coloring model and gas hard-core model, strong spatial mixing can be shown to hold only in certain ranges of the parameter space: in the Ising model, the temperature must be high enough; in the proper coloring model, the number of colours should be large enough; and in the gas hard-core model, the fugacity must be low enough. We hence see that the noisy voter model is extremely special, since it can be thought of as \emph{always} being in the high-temperature, highly disordered regime. 

\begin{remark}\label{infnvm}
 We also mention that the noisy voter model, when defined on the entire lattice $\bbZ^d$, is ergodic for \emph{all} $\delta > 0$, an easy consequence of Dobrushin's `$M < \epsilon$' criterion \cite{D}. This in particular means that the noisy voter model always has a unique infinite-volume invariant measure, again in contrast to other spin systems like the Ising model with zero magnetic field ($d \geq 2$): there exists a unique infinite-volume Gibbs measure only above the critical temperature.
\end{remark}

\section*{Acknowledgements}
This paper contains work done while writing my Master's thesis at the Ludwig-Maximilians-Universit\"at and Technische Universit\"at M\"unchen, supported by the Elitenetzwerk Bayern. I thank my supervisor, Nina Gantert (TU M\"unchen) for suggesting this topic, for the many enlightening discussions and for having taught me so much.

\end{document}